\documentclass{article}
\usepackage{graphicx}
\usepackage{eepic}
\usepackage{pst-all}
\usepackage{amsmath, amsthm,amsfonts, stmaryrd}
\usepackage[cyr]{aeguill}
\usepackage{listings}
\usepackage{psfrag}
\usepackage{epsfig}

\usepackage{cancel}

\renewcommand{\v}{\vec} 
\renewcommand{\c}{\times} 
\renewcommand{\o}{\circ}
\renewcommand{\b}{\mathbf} 

\newcommand{\w}{\vec{\omega}} 
\newcommand{\IM}{\vec{\textrm{Im}}} 
\newcommand{\RE}{\textrm{Re}} 

\renewcommand{\wp}{\v{\omega}'}

\newcommand{\bq}{\b{q}}

\newcommand{\half}{\frac{1}{2}}

\newcommand{\vbeg}{\left( \begin{array}{c} } 
\newcommand{\vend}{ \end{array} \right)} 

\definecolor{cBlue}{rgb}{0.0,0.0,1.0}
\definecolor{cViol}{rgb}{0.4,0.0,0.6}
\definecolor{cRed}{rgb}{0.7,0.0,0.3}
\definecolor{cGreen}{rgb}{0.3,0.5,0.0}
\definecolor{cGrey}{rgb}{0.2,0.2,0.2}
\definecolor{cGreyy}{rgb}{0.5,0.5,0.5}
\definecolor{cBlack}{rgb}{0.0,0.0,0.0}

\newcommand{\Blue}{\color{cBlue}} 
 
\newcommand{\Red}{\color{cRed}} 
\newcommand{\Green}{\color{cGreen}} 
 
\newcommand{\Greyy}{\color{cGreyy}} 
\newcommand{\Black}{\color{cBlack}}

\begin{document}

\title{Quaternions And Dynamics}
\author{Basile Graf \\ \texttt{basile.graf@epfl.ch}}
\date{February, 2007}
\maketitle
\thispagestyle{empty}

\subsection*{Abstract}
We give a simple and self contained introduction to quaternions and their practical usage in dynamics. The rigid body dynamics are presented in full details. In the appendix, some more exotic relations are given that allow to write more complex models, for instance, the one of a satellite with inertial wheels and expressed in a non-inertial reference frame. As it is well known, one nice advantage of quaternions over Euler angles, beside the usual arguments, is that it allows to write down quite complex dynamics completely by hand. 

\pagebreak

\tableofcontents

\pagebreak 

\graphicspath{{SubDocuments/Quaternions/}}

\renewcommand{\w}{\omega} 

\section{Quaternions}
\label{quat}

\subsection{Fundamentals}

Relation \eqref{eq1}, together with associativity and distributivity is all what we will use to derive the basic practical applications for quaternions.

\begin{equation}
\fbox{$i^2 = j^2 = k^2 = ijk = -1$}
\label{eq1}
\end{equation}

By left- and right-multiplication in the above equation, we can write

\begin{equation*}
\begin{tabular}{cc}
	\multicolumn{2}{c}{$i\ ijk = -jk = -i$}  \\
	\multicolumn{2}{c}{$ijk\ k = -ij = -k$}  \\
	$j \ jk = -k = ji$  &   $ij \ j = -i = kj$  \\
	$i \ ij = -j = ik$  &   $ji \ i = -j = -ki$ 
\end{tabular}
\end{equation*}

This shows the product is \emph{non commutative} and gives the basic multiplication rules:

\begin{equation}
\fbox{
\begin{tabular}{cc}
	$ij=k$ & $ji=-k$  \\
	$jk=i$ & $kj=-i$  \\
	$ki=j$ & $ik=-j$  \\
\end{tabular}}
\label{eq3}
\end{equation}

\subsection{Notations and Definitions}

A quaternion $q$ is a set of four parameters, a real value $q_0$ and three imaginary values $q_1i$, $q_2j$, $q_3k$ with $q_1, q_2, q_3 \in  \mathbb{R}$; it may be written
 
\begin{equation*}
q=q_0+q_1i+q_2j+q_3k .
\end{equation*}

However, this notation proves itself to be very unpractical. We will therefore use two different notations:

\begin{itemize}
	\item The quaternion $q$ as a pair of real value and vectorial imaginary value \\
			$q = (q_0, \vec{q})  \qquad\qquad   \RE\big\{ q \big\}=q_0  \qquad  \IM\big\{ q \big\}=\vec{q}=(q_1\ q_2\ q_3)^T$
	\item A column vector of four parameters \\
			$\b{q}=(q_0\ q_1\ q_2\ q_3)^T$
\end{itemize}

The conjugate $\bar{q}$ of $q$ is defined as

\begin{equation*}
\bar{q} = (q_0, -\vec{q})
\end{equation*}

and it's norm (a nonnegative real value) as

\begin{equation*}
|q| = |\b{q}| = \sqrt{q^2_0+q^2_1+q^2_2+q^2_3}.
\end{equation*}

The product of two quaternions written as pairs, as described in the next section will be noted with $\o$.

\subsection{Quaternion Product}

From the rules given in \eqref{eq3}, we may write the product of $q$ with $p$.

\begin{equation*}
(q_0 + q_1i +q_2j +q_3k) (p_0 + p_1i +p_2j +p_3k) \ =
\label{eq4}
\end{equation*}

\begin{equation*}
\begin{tabular}{ccccccccc}
    & \Greyy{$p_0 q_0$}    &$+$&   \Red{$q_0 p_1 \ i$}   &$+$&   \Red{$q_0 p_2 \ j$}   &$+$&   \Red{$q_0 p_3 \ k$}  & \\
$+$ & \Blue{$q_1 p_0 \ i$} &$+$&   \Green{$q_1 p_1 \ ii$}  &$+$&   $q_1 p_2 \ ij$  &$+$&   $q_1 p_3 \ ik$ & \\
$+$ & \Blue{$q_2 p_0 \ j$} &$+$&   $q_2 p_1 \ ji$  &$+$&   \Green{$q_2 p_2 \ jj$}  &$+$&   $q_2 p_3 \ jk$ & \\
$+$ & \Blue{$q_3 p_0 \ k$} &$+$&   $q_3 p_1 \ ki$  &$+$&   $q_3 p_2 \ kj$  &$+$&   \Green{$q_3 p_3 \ kk$} & $=$ \\
\end{tabular}
\label{eq5}
\end{equation*}

\begin{equation*}
\begin{tabular}{ccccccccc}
    & \Greyy{$p_0 q_0$}   &$-$&   \Green{$q_1 p_1$}   &$-$&   \Green{$q_2 p_2$}    &$-$&   \Green{$q_3 p_3$}   & \\
$+$ & \Blue{$(q_1 p_0$}   &$+$&   \Red{$q_0 p_1$}   &$+$&   $q_2 p_3$    &$-$&   $q_3 p_2)\ i$   & \\
$+$ & \Blue{$(q_2 p_0$}   &$+$&   \Red{$q_0 p_2$}   &$+$&   $q_3 p_1$    &$-$&   $q_1 p_3)\ j$   & \\
$+$ & \Blue{$(q_3 p_0$}   &$+$&   \Red{$q_0 p_3$}   &$+$&   $q_1 p_2$    &$-$&   $q_2 p_1)\ k$   & \\
\end{tabular}
\label{eq6}
\end{equation*}

\begin{equation}
q \circ p = (\Greyy{p_0 q_0} - \Green{\v{p}\cdot\v{q}} , \Red{q_0 \v{p}} + \Blue{p_0 \v{q}} + \Black{\v{q}\c\v{p}}).
\label{eq7}
\end{equation}

From \eqref{eq7} it turns out that

\begin{equation}
q \o \bar{q} = \bar{q} \o q = (|q|^2, \vec{0}) = |q|^2
\label{eqConj}
\end{equation}

and if $q$ is normed ($|q|=1$)

\begin{equation}
q \o \bar{q} = \bar{q} \o q = (1, \vec{0}) = Id.
\label{eqConjN}
\end{equation}

In \eqref{eq7} we also see that 

\begin{equation}
\overline{q \o p} = \bar{p} \o \bar{q}
\label{eqQoQbar}
\end{equation}

that is

\begin{equation*}
|q \o p|^2 = (q\o p) \o (\overline{q \o p}) = q\o \underbrace{p \o \bar{p}}_{|p|^2} \o \bar{q} = |p|^2 (q\o \bar{q}) = |q|^2|p|^2
\end{equation*}
\begin{equation}
|q \o p|  = |q||p|.
\label{eqNormProd}
\end{equation}

\subsection{Quaternions and Spatial Rotations}

First, note the following relations

\begin{gather*}
(\vec{u} \c \vec{v}) \c \vec{w} = (\vec{u} \cdot \vec{w}) \vec{v}   -   (\vec{v} \cdot \vec{w}) \vec{u} \\
\sin^2\frac{\varphi}{2} = \frac{1-\cos\varphi}{2}  \qquad    \cos^2\frac{\varphi}{2} = \frac{1+\cos\varphi}{2}.
\end{gather*}

From now on, $q$ will generally represent a normed quaternion ($|q|=1$) involved in a rotation. Let's now place a vector $\vec{x}\in \mathbb{R}^3$ in the imaginary part of a quaternion $x$ and see what happens with it in the following relation   

\begin{equation*}
x' = \bar{q} \o x \o q  \qquad \qquad   x = (0,\vec{x}) \qquad  q = (q_0, \vec{q}) .  
\end{equation*}

Using \eqref{eq7}

\begin{gather*}
x' = (\vec{q} \cdot \vec{x},\  q_0 \vec{x} - \vec{q} \c \vec{x}) \o q \\
= (  \underbrace{(\vec{q} \cdot \vec{x})q_0 - (q_0 \vec{x} - \vec{q} \c \vec{x}) \cdot \vec{q}}_{\RE\{ x'  \}},\  
\underbrace{(\vec{q} \cdot \vec{x}) \vec{q} + q_0(q_0 \vec{x} - \vec{q} \c \vec{x}) +(q_0 \vec{x} - \vec{q} \c \vec{x}) \c \vec{q}}_{\IM\{  x' \}}  )
\end{gather*}

\begin{align*}
\RE\big\{ x' \big\} &= (\vec{q} \cdot \vec{x})q_0 - q_0(\vec{x} \cdot \vec{q}) - (\vec{q} \c \vec{x}) \cdot \vec{q} = 0 \\
                    &\Rightarrow \ x' = (0,\vec{x}'), \\
\IM\big\{ x' \big\} &= \vec{x}' \\
										&= (\vec{q} \cdot \vec{x}) \vec{q} + q_0^2\vec{x} - q_0(\vec{q}\c\vec{x}) + q_0(\vec{x}\c\vec{q}) - (\vec{q}\c\vec{x}) \c\vec{q} \\
                    &= (\vec{q} \cdot \vec{x}) \vec{q} + q_0^2\vec{x} + 2 q_0(\vec{x}\c\vec{q}) - (\vec{q}\c\vec{x}) \c\vec{q} \\
    &= (\vec{q} \cdot \vec{x}) \vec{q} + q_0^2\vec{x} + 2 q_0(\vec{x}\c\vec{q}) - (\vec{q}\cdot\vec{q})\vec{x} + (\vec{x}\cdot\vec{q})\vec{q} \\
    &= 2(\vec{q} \cdot \vec{x}) \vec{q} + q_0^2\vec{x} + 2 q_0(\vec{x}\c\vec{q}) - (\vec{q}\cdot\vec{q})\vec{x}.
\end{align*}

A valid normed quaternion ($|q|=\sqrt{(q_0^2+q_1^2+q_2^2+q_3^2)}=1$) would be

\begin{equation*}
q = (q_0, \vec{q}) = (\cos\frac{\varphi}{2},\ \sin\frac{\varphi}{2}\vec{n}) \qquad \qquad |\vec{n}|=1.
\end{equation*}

In this case, $\vec{x}'$ becomes 

\begin{align*}
\vec{x}'  &=  
2 \sin^2\frac{\varphi}{2} (\vec{n}\cdot\vec{x})\vec{n} + \cos^2\frac{\varphi}{2} \vec{x} + 2 \cos\frac{\varphi}{2} \sin\frac{\varphi}{2} (\vec{x}\c\vec{n}) - \sin^2\frac{\varphi}{2} \vec{x} \\
         &= (1-\cos\varphi)(\vec{n}\cdot\vec{x}) \vec{n}  + \cos\varphi\ \vec{x} + \sin\varphi\ (\vec{x}\c\vec{n}).
\end{align*}

This last relation is the formula for a rotation by an angle $\varphi$ around a normed axis vector $\vec{n}$, as can be shown with the following figure as follows:

\begin{center}
\psfrag{0}[l][c][1][0]{$0$}
\psfrag{n}[l][c][1][0]{$\vec{n}$}
\psfrag{xnn}[l][c][1][29]{$(\vec{x}\cdot\vec{n})\vec{n}$}
\psfrag{x}[l][c][1][0]{$\vec{x}$}
\psfrag{xp}[l][c][1][0]{$\vec{x}'$}
\psfrag{v1}[l][c][1][0]{$\vec{v}_1$}
\psfrag{v2}[l][c][1][0]{$\vec{v}_2$}
\psfrag{v3}[l][c][1][0]{$\vec{v}_3$}
\psfrag{phi}[l][c][1][0]{$\varphi$}
\epsfig{file=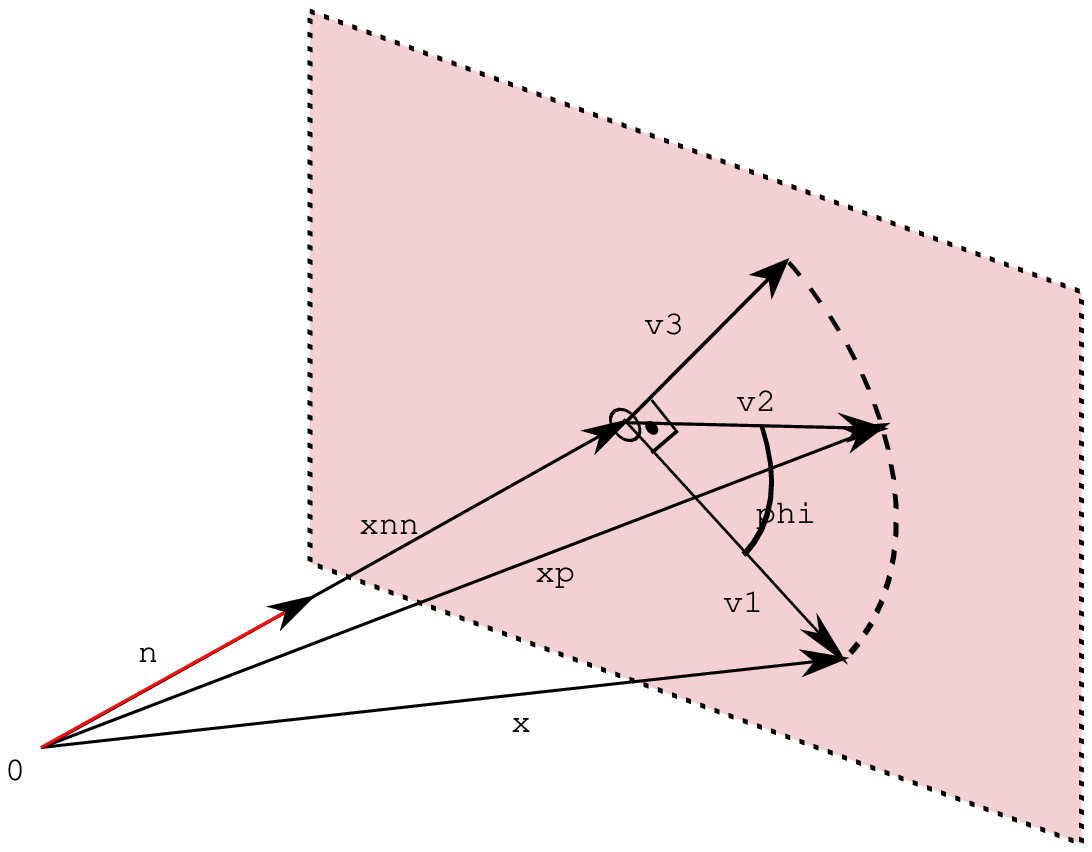, width=8cm}
\end{center}

\begin{align*}
\vec{v}_2 &= \cos\varphi\ \vec{v}_1 + \sin\varphi\ \vec{v}_3 \\
\vec{v}_1 &= \vec{x} - (\vec{x}\cdot\vec{n})\vec{n} \\
\vec{v}_3 &= \vec{v}_1 \c \vec{n}  \\
          &= (\vec{x} - (\vec{x}\cdot\vec{n})\vec{n}) \c \vec{n} \\
          &= (\vec{x}\c\vec{n}) - (\vec{x}\cdot\vec{n})\underbrace{\vec{n}\c\vec{n}}_{\vec{0}} \\
\Rightarrow\ \vec{v}_2 &= \cos\varphi\ (\vec{x} - (\vec{x}\cdot\vec{n})\vec{n}) + \sin\varphi\ (\vec{x}\c\vec{n}) \phantom{\underbrace{1}_{\vec{1}}}\\
\vec{x}'  &= (\vec{x}\cdot\vec{n})\vec{n} + \vec{v}_2 \\
          &= (\vec{x}\cdot\vec{n})\vec{n} + \cos\varphi\ (\vec{x} - (\vec{x}\cdot\vec{n})\vec{n}) + \sin\varphi\ (\vec{x}\c\vec{n}) \\
          &= (1-\cos\varphi)(\vec{n}\cdot\vec{x}) \vec{n}  + \cos\varphi\ \vec{x} + \sin\varphi\ (\vec{x}\c\vec{n}).
\end{align*}

Moreover

\begin{gather*}
x' = \bar{q} \o x \o q \\
q \o x' \o \bar{q} = \underbrace{q\o\bar{q}}_{(1,\vec{0})} \o x \o \underbrace{q\o\bar{q}}_{(1,\vec{0})}.
\end{gather*}

Thus we have the relations for the rotation and its inverse

\begin{equation}
\fbox{$x' = \bar{q} \o x \o q \\   \qquad  x = q \o x' \o \bar{q}$}.
\end{equation}

\subsection{Quaternions and Rotation Velocity}

We will now derive the relation between the rotational velocity vector and the quaternion time derivative. $\vec{x}'$ is any constant vector within the \emph{body (rotating) reference frame} and $\vec{x}$ is the same vector in the \emph{fixed reference frame}. As seen before, both vectors can be put in relation with

\begin{equation*}
\begin{tabular}{ccc}
$x = q \o x' \o \bar{q}$ &  & $x' = \bar{q} \o x \o q$.
\end{tabular}
\label{eq8}
\end{equation*}

Applying the time derivative to $x=(0,\vec{x})$, with $x'=(0,\vec{x}')$ and $\dot{\vec{x}}'=\v{0}$, we get

\begin{gather*}
\dot{x} = \dot{q} \o x' \o \bar{q}  +  q \o x' \o \dot{\bar{q}}  \\
\dot{x} = \dot{q} \o \bar{q} \o x \o \underbrace{q \o \bar{q}}_{Id}  +  \underbrace{q \o \bar{q}}_{Id} \o x \o q \o \dot{\bar{q}}
\end{gather*}
\begin{equation}
\dot{x} = \dot{q} \o \bar{q} \o x  +   x \o q \o \dot{\bar{q}}  
\label{eq9.3}
\end{equation}

and from \eqref{eq7}

\begin{gather*}
\dot{q} \o \bar{q} = ( \underbrace{\dot{q_0}q_0 + \dot{\v{q}}\cdot\v{q}}_{\varoast}, -\dot{q_0}\v{q} + q_0\dot{\v{q}} - \dot{\v{q}} \c \v{q}) \\
\varoast = q_0\dot{q_0} + q_1\dot{q_1} + q_2\dot{q_2} + q_3\dot{q_3} = \b{q} \cdot \dot{\b{q}} = 0
\end{gather*}

because $|\b{q}|=1$. That is

\begin{equation}
\dot{q} \o \bar{q} = (0,\v{\nu})  \qquad \textrm{and similarly} \qquad \bar{q} \o \dot{q} = (0,\v{-\nu}). 
\label{eq10}
\end{equation}

\subsubsection{Rotation Velocity in Fixed Reference Frame $\w$}
From \eqref{eq9.3} and \eqref{eq10} and using \eqref{eq7} we have

\begin{align*}
\dot{x} & = (0,\v{\nu}) \o x - x \o (0,\v{-\nu})  \\
\dot{\vec{x}} & = \vec{\nu} \c \vec{x} - \vec{x} \c \vec{\nu} = 2 \vec{\nu} \c \vec{x}
\end{align*}
and from \eqref{eqNormProd} 
\begin{equation*}
|\dot{\vec{x}}|=|2\vec{\nu}||\vec{x}|  \qquad \Rightarrow \qquad    \vec{\nu} \bot \vec{x}
\end{equation*}

If $\vec{x}$ undergoes a pure rotation, we know that

\begin{equation*}
\dot{\vec{x}} = \vec{\w} \c \vec{x}    \qquad \textrm{and} \qquad   \vec{\w} \bot \vec{x}
\end{equation*}

thus

\begin{equation}
\fbox{ $\w = (0,\vec{\w}) = 2(0,\vec{\nu}) = 2 \dot{q} \o \bar{q} $}.
\end{equation}

And right-multiplication by $q$

\begin{equation*}
\w \o q= 2 \dot{q} \o \underbrace{\bar{q} \o q}_{Id} \ \  \Rightarrow \ \  \w \o q = 2 \dot{q}
\end{equation*}

\begin{equation}
\fbox{$ \dot{q} = \frac{1}{2} \w \o q $}.
\end{equation}

\subsubsection{Rotation Velocity in Body Reference Frame $\w'$}

\begin{align*}
             \w'      & = \bar{q} \o \w \o q    \qquad \textrm{with} \qquad    \w = 2 \dot{q} \o \bar{q} \\
\Rightarrow  \ \ \w'  & = 2 \bar{q} \o \dot{q} \o \underbrace{\bar{q} \o q}_{Id}
\end{align*}

\begin{equation}
\fbox{ $\w' = 2 \bar{q} \o \dot{q} $}.
\end{equation}

And left-multiplication by $q$

\begin{equation*}
q \o \w' = 2 \underbrace{q \o \bar{q}}_{Id} \o \dot{q} = 2 \dot{q}
\end{equation*}

\begin{equation}
\fbox{$ \dot{q} = \frac{1}{2} q \o \w' $}.
\end{equation}

\subsubsection{Matrix-Product Notation for $\w$}
\label{secMat1}

From

\begin{equation*}
\w = 2 \dot{q} \o \bar{q} 
\end{equation*}

and using \eqref{eq7}

\begin{align*}
\vec{\w} &= \IM\big\{{2 \dot{q} \o \bar{q}} \big\} = 2(\Red{-\dot{q_0}\vec{q}} \Blue{ + q_0\dot{\vec{q}}} \Black{- \dot{\vec{q}} \c \vec{q})} \\
         &= 2 \underbrace{\left( \begin{array}{cccc}
\Red{-q_1} &  \Blue{q_0} & -q_3 & q_2 \\
\Red{-q_2} &  q_3 & \Blue{q_0}  & -q_1 \\
\Red{-q_3} & -q_2 & q_1  & \Blue{q_0} \\
\end{array} \right)}_{E} \left( \begin{array}{c}
\dot{q_0} \\
\dot{q_1} \\
\dot{q_2} \\
\dot{q_3} \\
\end{array} \right)
\end{align*}

\begin{equation*}
\vec{\w} = 2E\dot{\b{q}} .
\end{equation*}

Changing the sign and inverting the cross product allows to make an other identification

\begin{equation*}
\vec{\w} = -2(-q_0\dot{\vec{q}} + \dot{q_0}\vec{q} - \vec{q} \c \dot{\vec{q}}) 
\end{equation*}

\begin{equation*}
\vec{\w} = -2\dot{E}\b{q} .
\end{equation*}

So the rotation velocity vector in the \emph{fixed reference} frame can be written as

\begin{equation}
\fbox{$ \vec{\w} = 2E\dot{\b{q}} = -2\dot{E}\b{q} $} .
\end{equation}

And from

\begin{equation*}
\dot{q} = \frac{1}{2} \w \o q     \qquad \qquad  \w = (0,\vec{\w}) \ \Rightarrow \ \w_0=0
\end{equation*}

one can similarly find

\begin{equation*}
\dot{\b{q}} = \frac{1}{2} \left( \begin{array}{c}
(-\vec{\w} \cdot \vec{q})\\
(q_0\vec{\w} + \vec{\w} \c \vec{q})
\end{array} \right) = \frac{1}{2} E^T\vec{\w}
\end{equation*}

\begin{equation}
\fbox{$ \dot{\b{q}} = \frac{1}{2} E^T\vec{\w} $} .
\label{eq_dotq_w}
\end{equation}

\subsubsection{Matrix-Product Notation for $\w'$}
\label{secMat2}

From

\begin{equation*}
\w' = 2 \bar{q} \o \dot{q} 
\end{equation*}

and using \eqref{eq7}

\begin{align*}
\vec{\w'} &= \IM\big\{{2 \bar{q} \o \dot{q}} \big\} = 2(\Blue{q_0\dot{\vec{q}}} \Red{ - \dot{q_0}\vec{q}} \Black{- \vec{q} \c \dot{\vec{q}})} \\
         &= 2 \underbrace{\left( \begin{array}{cccc}
\Red{-q_1} &  \Blue{q_0} & q_3 & -q_2 \\
\Red{-q_2} & -q_3 & \Blue{q_0}  & q_1 \\
\Red{-q_3} & q_2 & -q_1  & \Blue{q_0} \\
\end{array} \right)}_{G} \left( \begin{array}{c}
\dot{q_0} \\
\dot{q_1} \\
\dot{q_2} \\
\dot{q_3} \\
\end{array} \right)
\end{align*}

\begin{equation*}
\vec{\w'} = 2G\dot{\b{q}} .
\end{equation*}

Changing the sign and inverting the cross product allows to make an other identification

\begin{equation*}
\vec{\w'} = -2(\dot{q_0}\vec{q} - q_0\dot{\vec{q}} - \dot{\vec{q}} \c \vec{q}) 
\end{equation*}

\begin{equation*}
\vec{\w'} = -2\dot{G}\b{q} .
\end{equation*}

So the rotation velocity vector in the \emph{body reference} frame can be written as

\begin{equation}
\fbox{$ \vec{\w'} = 2G\dot{\b{q}} = -2\dot{G}\b{q} $} .
\label{eqGq}
\end{equation}

And from

\begin{equation*}
\dot{q} = \frac{1}{2} q \o \w'     \qquad \qquad  \w' = (0,\vec{\w'}) \ \Rightarrow \ \w'_0=0
\end{equation*}

one can similarly find

\begin{equation*}
\dot{\b{q}} = \frac{1}{2} \left( \begin{array}{c}
(-\vec{q} \cdot \vec{\w'})\\
(q_0\vec{\w'} + \vec{q} \c \vec{\w'})
\end{array} \right) = \frac{1}{2} G^T\vec{\w'}
\end{equation*}

\begin{equation}
\fbox{$ \dot{\b{q}} = \frac{1}{2} G^T\vec{\w'} $} .
\label{eqqdot}
\end{equation}

\subsubsection{Rotation Matrix $R$}

We already have

\begin{minipage}[left]{0.49 \textwidth}
	\begin{gather*}
	   	 \vec{\w} = 2E\dot{\b{q}} = -2\dot{E}\b{q}  \\
		\dot{\b{q}} = \frac{1}{2} E^T\vec{\w}
	\end{gather*}
\end{minipage}
\begin{minipage}[right]{0.49 \textwidth}
	\begin{gather*}
	   	 \vec{\w'} = 2G\dot{\b{q}} = -2\dot{G}\b{q}  \\
		\dot{\b{q}} = \frac{1}{2} G^T\vec{\w'}
	\end{gather*}
\end{minipage}

So we can write

\begin{minipage}[left]{0.49 \textwidth}
	\begin{align*}
		\vec{\w} &= 2E\dot{\b{q}}   \\
		         &= 2E(\frac{1}{2} E^T\vec{\w})   \\
		         &= EE^T \vec{\w}   \\
		&\Rightarrow \ \fbox{$EE^T = Id$} .
	\end{align*}
\end{minipage}
\begin{minipage}[right]{0.49 \textwidth}
	\begin{align*}
		\vec{\w'} &= 2G\dot{\b{q}}   \\
		         &= 2G(\frac{1}{2} G^T\vec{\w'})   \\
		         &= GG^T \vec{\w'}   \\
		&\Rightarrow \ \fbox{$GG^T = Id$} .
	\end{align*}
\end{minipage} \\

And by mixing both sides

\begin{gather*}
	\vec{\w'} = 2G\dot{\b{q}} = 2G(\frac{1}{2} E^T\vec{\w}) = GE^T \vec{\w} \\
	\vec{\w} = 2E\dot{\b{q}} = 2E(\frac{1}{2} G^T\vec{\w'}) = EG^T \vec{\w'} .
\end{gather*}

We shall now remember that $\v{\w}$ is a vector in the \emph{fixed reference frame} and that $\v{\w}'$ is the same vector in the \emph{body reference frame}, that is $\vec{\w} = R \vec{\w'}$. By comparing with the previous two results, we find

\begin{equation}
\fbox{$ R = EG^T $} \qquad \textrm{and} \qquad \fbox{$R^{-1}=R^T=GE^T $} .
\label{eqRot}
\end{equation}

\subsubsection{$E\b{p}$ and $G\b{p}$}

From the identifications made in sections \ref{secMat1} and \ref{secMat2}, we can see that the general meaning the product of $E$ and $G$ with any quaternion $\b{p}$ is

\begin{equation}
E\b{p} = \IM\big\{{p \o \bar{q}}\big\}    \qquad \qquad    G\b{p} = \IM\big\{{\bar{q} \o p} \big\} .
\label{eqEpGp}
\end{equation}

And from

\begin{equation*}
q \o \bar{q} = \bar{q} \o q = (|q|, \vec{0}) = (1,\vec{0})    
\end{equation*}

it follows

\begin{equation*}
\fbox{$E\b{q} = \vec{0}$}    \qquad \qquad    \fbox{$G\b{q} = \vec{0}$} .
\end{equation*}

\subsubsection{One Last Relation}

For any $\vec{v}$ and due to associativity

\begin{align*}
\underbrace{(0,\vec{\w'})}_{2 \bar{q} \o \dot{q}} \o (0,\vec{v})  &=  (-\vec{\w'}\cdot\vec{v} , \vec{\w'}\c \vec{v}) \\
           &= 2 \bar{q} \o \dot{q} \o v \\
           &=  2(\Greyy{q_0\dot{q_0}} \Green{+\vec{q} \cdot \dot{\vec{q}}},  \Blue{q_0\dot{\vec{q}}} \Red{ - \dot{q_0}\vec{q}} \Black{- \vec{q} \c \dot{\vec{q}})} \o v =
2 \bar{q} \o (\Greyy{\dot{q_0}v_0} \Green{-\dot{\vec{q}} \cdot \vec{v}},  \Blue{\dot{q_0}\vec{v}} \Red{ + v_0\dot{\vec{q}}} 
\Black{+ \dot{\vec{q}} \c \vec{v}}) \\
         &\equiv 2 \left( \begin{array}{cccc}
\Greyy{q_0} & \Green{q_1} & \Green{q_2} & \Green{q_3} \\
\hline
\Red{-q_1} &  \Blue{q_0} & q_3 & -q_2 \\
\Red{-q_2} & -q_3 & \Blue{q_0}  & q_1 \\
\Red{-q_3} & q_2 & -q_1  & \Blue{q_0} \\
\end{array} \right) 
\left( \begin{array}{c|cccc}
\Greyy{\dot{ q_0}} &  \Green{\dot{-q_1}} & \Green{\dot{-q_2}} & \Green{\dot{-q_3}}\\
\Red{\dot{ q_1}} &  \Blue{\dot{ q_0}} & \dot{-q_3} & \dot{ q_2}\\
\Red{\dot{ q_2}} &  \dot{ q_3} & \Blue{\dot{ q_0}} & \dot{-q_1}\\
\Red{\dot{ q_3}} &  \dot{-q_2} & \dot{ q_1} & \Blue{\dot{ q_0}}\\
\end{array} \right) 
\left( \begin{array}{c}
0 \\
\hline
v_1 \\
v_2 \\
v_3 \\
\end{array} \right) \\
         &= 2 \left( \begin{array}{c}
         \b{q}^T \\
         G  \\
         \end{array} \right) 
         \bigg(\begin{array}{cc}
         \dot{\b{q}} & \dot{G}^T \\
         \end{array} \bigg) 
         \bigg(\begin{array}{c}
         0 \\
         \vec{v} \\
         \end{array} \bigg) = 
         \bigg(\begin{array}{c}
         -\vec{\w'}\cdot\vec{v} \\
         \vec{\w'}\c \vec{v}\\
         \end{array} \bigg)\\
     &\Rightarrow \ 2G\dot{G}^T\vec{v}=\Omega' \vec{v} = \vec{\w'} \c \vec{v} .
\end{align*}

Comparing with \eqref{eqGq}, we conclude that

\begin{equation}
\fbox{$ \Omega'= 2G\dot{G}^T = -2\dot{G}G^T \qquad \textrm{and} \qquad  \Omega' \vec{v} = \vec{\w'} \c \vec{v} $} .
\label{eqCross}
\end{equation}

\subsubsection{Relations Summary}

The following table summaries the developed relations. $q$ is always a normed quaternion, that is $q_0^2+q_1^2+q_2^2+q_3^2=1$. \\
\begin{center}
\begin{tabular}{| p{0.249\textwidth} | p{0.249\textwidth} | p{0.249\textwidth} | p{0.249\textwidth} |}
\hline
\multicolumn{2}{|c|}{Quaternion notation}&
\multicolumn{2}{|c|}{Matrix notation}\\
\hline
\multicolumn{1}{|c|}{\emph{Fixed ref}} &
\multicolumn{1}{|c|}{\emph{Body ref}} &
\multicolumn{1}{|c|}{\emph{Fixed ref}} &
\multicolumn{1}{|c|}{\emph{Body ref}} \\
\hline
$x = q \o x' \o \bar{q}     \phantom{\Bigg(}$  &
$x' = \bar{q} \o x \o q$  &
\begin{tabular}{l}
	$\vec{x}=R\vec{x}'$ \\
	$R=EG^T$
\end{tabular}&
\begin{tabular}{l}
	$\vec{x}'=R^T\vec{x}$  \\
	$R^T=R^{-1}=GE^T$ 
\end{tabular}\\
\hline
$\w=(0,\vec{\w})=2\dot{q}\o\bar{q}     \phantom{\Bigg(}$  &
$\w'=(0,\vec{\w}')=2\bar{q}\o\dot{q}$  &
$\vec{\w}=2E\dot{\b{q}}=-2\dot{E}\b{q}$  &
$\vec{\w}'=2G\dot{\b{q}}=-2\dot{G}\b{q}$ \\
\hline
$\dot{q}=\frac{1}{2}\w \o q     \phantom{\Bigg(}$  &
$\dot{q}=\frac{1}{2} q \o \w'$  &
$\dot{\b{q}}=\frac{1}{2}E^T\vec{\w}$  &
$\dot{\b{q}}=\frac{1}{2}G^T\vec{\w}'$  \\
\hline
&
&
$EE^T=Id     \phantom{\Bigg(}$  &
$GG^T=Id$  \\
\hline
\multicolumn{2}{|c|}{$q\o\bar{q} = \bar{q}\o q = (|q|,\vec{0})     \phantom{\Bigg(}$} &
$E\b{q}=\vec{0}$ &
$G\b{q}=\vec{0}$ \\
\hline
$     \phantom{\Bigg(}$   &
\begin{tabular}{l}
	$(0,\vec{\w}) \o (0,\vec{v}) = $ \\
	$(-\vec{\w'}\cdot\vec{v},\vec{\w'}\c\vec{v})$ 
\end{tabular}&
&
\begin{tabular}{l}
	$\Omega' = 2G\dot{G}^T$ \\
	$\ \ \ \  = -2\dot{G}G^T$  \\
	$\Omega'\vec{v} = \vec{\w'}\c\vec{v}$ 
\end{tabular}\\
\hline
\end{tabular}
\end{center}

\begin{equation*}
E =
\left( \begin{array}{cccc}
-q_1 &  q_0 & -q_3 & q_2  \\
-q_2 &  q_3 & q_0  & -q_1 \\
-q_3 & -q_2 & q_1  & q_0  \\
\end{array} \right)
\qquad \qquad
G = 
\left( \begin{array}{cccc}
-q_1 &  q_0 & q_3 & -q_2 \\
-q_2 & -q_3 & q_0  & q_1 \\
-q_3 & q_2 & -q_1  & q_0 \\
\end{array} \right)
\end{equation*}

\subsection{Rigid Body Rotational Dynamics}
\label{quat_rigid_dyn_IN}

We now will have a look at the dynamics of a freely rotating rigid body to which a momentum $\v T'$ is applied. Translation of the body will not be discussed (it can be decoupled from the dynamics of rotation and is fairly easy). We will also consider a potential free system, so that the Lagrangian resumes to the rotational kinetic energy only

\begin{equation}
L=E_{\textrm{rot}}=\frac{1}{2}\vec{\w}'^TJ\vec{\w}'.
\end{equation}

Using the quaternion $\b{q}$ as coordinates and with the constraint $C=\b{q}^T\b{q}=1$, Lagrangian dynamics gives

\begin{equation}
\frac{d}{dt}\frac{\partial L}{\partial \dot{\b{q}}} - \frac{\partial L}{\partial \b{q}} = \b{F_q} + {\lambda} \frac{\partial C}{\partial \b{q}}.
\label{eqLDyn}
\end{equation}

$\b{F_q}$ is the 4-vector of generalized forces which will be expressed in term of applied torque later. $\b{\lambda}$ is the Lagrangian multiplier used to satisfy the constraint $C$.

\subsubsection{Derivatives of $L$}

Note the following reminder

\begin{center}
\begin{tabular}{c}
$\frac{\displaystyle \partial A\b{x}}{\displaystyle \partial \b{x}} = A      \phantom{\Bigg(}$ \\
$\frac{\displaystyle \partial \b{a}^T\b{x}}{\displaystyle \partial \b{x}} = \frac{\displaystyle \partial \b{x}^T\b{a}}{\displaystyle \partial \b{x}} = \b{a}$ \\
$\frac{\displaystyle \partial \b{x}^T A\b{x}}{\displaystyle \partial \b{x}} = (A^T +A)\b{x} \stackrel{\textrm{if }A=A^T}{=} 2A\b{x}     \phantom{\Big(}$ \\
(written as column vectors)\\
$(AB)^T=B^TA^T \ .  \phantom{\Bigg(}$\\
\end{tabular}
\end{center}

We will now derive each term of the left side of \eqref{eqLDyn}. First, let us rewrite $L$ in two different ways

\begin{equation*}
L=\frac{1}{2}\vec{\w}'^TJ\vec{\w}' = 2(G\dot{\b{q}})^T J (G\dot{\b{q}}) = 2 (\dot{G}\b{q})^T J (\dot{G}\b{q})
\end{equation*}

and grouping around $J$

\begin{equation*}
L=\frac{1}{2}\vec{\w}'^TJ\vec{\w}' = 2\dot{\b{q}}^T (G^T J G) \dot{\b{q}} = 2 \b{q}^T (\dot{G}^T J \dot{G}) \b{q}.
\end{equation*}

Because $J$ is symmetric, $(G^T J G)$ and $(\dot{G}^T J \dot{G})$ are also symmetric. So we have

\begin{equation}
\frac{\partial L}{\partial \b{q}} = 4 \dot{G}^T J \dot{G} \b{q} = 2 \dot{G}^T J \underbrace{(2\dot{G}\b{q})}_{-\vec{\w}'} =
-2 \dot{G}^T J \vec{\w}',
\end{equation}

\begin{equation*}
\frac{\partial L}{\partial \dot{\b{q}}} = 4 G^T J G \dot{\b{q}} = 2 G^T J \underbrace{(2G\dot{\b{q}})}_{\vec{\w}'} =
2 G^T J \vec{\w}'
\end{equation*}

and

\begin{equation}
\frac{d}{dt} \frac{\partial L}{\partial \dot{\b{q}}} = \frac{d}{dt} (2 G^T J \vec{\w}') = 2 \dot{G}^T J \vec{\w}' + 2 G^T J \dot{\vec{\w}}'.
\end{equation}

\subsubsection{Generalized Forces}

A way to find the generalized force $\b{F_c}$ relative to the coordinates $\b{c}$ is to identify it in 

\begin{equation*}
\delta W = \b{F_c} \cdot \delta \b{c} .
\end{equation*}

(A simple example is the case of a pure translation $\delta \vec{x}$ of a particle, on which a force $\vec{F}$ is applied. The work is then $\delta W = \b{F}_{\vec{x}} \cdot \delta \vec{x} = \vec{F} \cdot \delta \vec{x}$. So the generalized force $\b{F}_{\vec{x}}$ is simply $\vec{F}$ in this case.) \\

For a rotation of a rigid body by an angle $\delta \varphi$ around an axis $\vec{n}$ with an applied torque $\vec{T}'$, the work can be written as

\begin{equation}
\delta W = (\vec{n} \cdot \vec{T}') \delta \varphi  \qquad \qquad   |\vec{n}|=1.
\label{eqWTfi}
\end{equation}

This small attitude change can be represented on one side as a small variation $\delta q$ of the coordinate quaternion $q$ and, on the other side, as a rotation quaternion $q_\delta$ operating from the current attitude represented by $q$ (i.e. a composition). That is

\begin{gather*}
q + \delta q = q \o q_\delta    \\
|q|=1 \qquad   |q_\delta|=1  \qquad    |\delta q| \ll 1.
\end{gather*}

We do not need to consider the fact that the variation $\delta q$ has to preserve the norm of $q$, because it will automatically be satisfied by introducing a constraint in the Lagrange formulation.\\

On one side we can write

\begin{gather*}
q + \delta q = q \o q_\delta    \\
\underbrace{\bar{q} \o q}_{(1,\vec{0})} + \bar{q} \o \delta q = \underbrace{\underbrace{\bar{q} \o q}_{(1,\vec{0})} \o q_\delta}_{q_\delta}    
\end{gather*}

\begin{equation}
\Rightarrow \ q_\delta = (1,\vec{0}) + \bar{q} \o \delta q.
\label{eq_qdelta}
\end{equation}

On the other side

\begin{equation*}
q_\delta = (\cos{\frac{\delta \varphi}{2}},\  \sin{\frac{\delta \varphi}{2}}\ \vec{n}).
\end{equation*}

Looking at the imaginary part

\begin{gather*}
\IM\big\{ q_\delta \big\} = \IM\big\{ \bar{q} \o \delta q \big\} = \sin{\frac{\delta \varphi}{2}}\ \vec{n} \approx  \frac{\delta \varphi}{2}\ \vec{n}
\end{gather*}

comparing with \eqref{eqWTfi}

\begin{equation*}
\Rightarrow \ \delta W = 2 \  \IM\big\{ \bar{q} \o \delta q \big\} \cdot \vec{T}'
\end{equation*}

and from \eqref{eqEpGp}

\begin{gather*}
\IM\big\{ \bar{q} \o \delta q \big\} = G \b{\delta q} \\
\Rightarrow \ \delta W = 2(G \b{\delta q}) \cdot \vec{T}' = 2 \vec{T}'^T (G \b{\delta q})  = 2(G^T \vec{T}')^T \b{\delta q} = \underbrace{2(G^T \vec{T}')}_{\b{F_q}} \cdot \b{\delta q}
\end{gather*}

\begin{equation}
\Rightarrow \ \fbox{$ \b{F_q} = 2G^T \vec{T}'$}.
\label{eqForce}
\end{equation}

\subsubsection{Dynamics}
\label{quat_rigid_dyn_IN_3}

We have now everything to write the dynamics

\begin{gather*}
\frac{d}{dt}\frac{\partial L}{\partial \dot{\b{q}}} - \frac{\partial L}{\partial \b{q}} = \b{F_q} + \lambda \frac{\partial C}{\partial \b{q}} \\
4 \dot{G}^T J \vec{\w}' + 2 G^T J \dot{\vec{\w}}' = 2G^T \vec{T}' + \lambda \b{q}. \\
\end{gather*}

Left-multiplying by $G$

\begin{gather*}
\underbrace{4 G \dot{G}^T}_{2\Omega'} J \vec{\w}' + 2 \underbrace{GG^T}_{Id} J \dot{\vec{\w}}' = 2 \underbrace{GG^T}_{Id} \vec{T}' + \lambda \underbrace{G \b{q}}_{\vec{0}} \\
\Omega' J \vec{\w}' +  J \dot{\vec{\w}}' =  \vec{T}' \\
\vec{\w}' \c J \vec{\w}' +  J \dot{\vec{\w}}' =  \vec{T}' \\
J \dot{\vec{\w}}' =  \vec{T}' - \vec{\w}' \c J \vec{\w}'.  
\end{gather*}

This last relation is nothing else than the Euler equation of motion for rotating body. Together with \eqref{eqqdot} we obtain the complete dynamics

\begin{equation}
\begin{tabular}{|rcl|}
	\hline
	$\dot{\vec{\w}}'$ & $=$ & $J^{-1} \vec{T}' - J^{-1} (\vec{\w}' \c J \vec{\w}')   \phantom{\Big(}$ \\
	$\dot{\b{q}}$     & $=$ & $\frac{1}{2} G^T\vec{\w}'.   \phantom{\Big(}$ \\
	\hline
\end{tabular}
\label{rigid_euler_IN}
\end{equation}


%

%
%
%
%
%
%


%

%

%


\appendix
\pagebreak 

\graphicspath{{SubDocuments/QuadraticFormDerivative/}}

\section{Derivatives and Quaternions}

\subsection{Quadratic Form Derivative by a Quaternion}
\label{dQuadForm_dq}

In order to be able to derive the Lagrangian by the components of $\bq$ for a non-inertial quaternion model, one needs to perform things like

\begin{equation*}
\frac{\partial (\v{v}^TR\v{w})}{\partial \bq} \ \ \ ,  \qquad \qquad   \frac{\partial (\v{v}^TR^T\v{w})}{\partial \bq}
\end{equation*}

and also

\begin{equation*}
\frac{\partial (\v{u}^TRJR^T\v{u})}{\partial \bq}.
\end{equation*}

But because $R=EG^T$ and 

\begin{equation*}
E =
\left( \begin{array}{cccc}
-q_1 &  q_0 & -q_3 & q_2  \\
-q_2 &  q_3 & q_0  & -q_1 \\
-q_3 & -q_2 & q_1  & q_0  \\
\end{array} \right)
\qquad \qquad
G = 
\left( \begin{array}{cccc}
-q_1 &  q_0 & q_3 & -q_2 \\
-q_2 & -q_3 & q_0  & q_1 \\
-q_3 & q_2 & -q_1  & q_0 \\
\end{array} \right)
\end{equation*}

the matrix of the quadratic form to be derived is not constant in $\bq$. This implies that these operations are no more trivial. However, thanks to the particular form of the dependance of $R$ in the components of $\bq$, higher order tensors can be avoided, as shown in the following.

\subsubsection{"Single $R$" Quadratic Form}

By computing the quadratic form and taking the partial derivatives we get (placing them in a column vector)

\begin{equation*}
\frac{\partial (\v{v}^TR\v{w})}{\partial \bq} =\left( \frac{\partial (\v{v}^TR\v{w})}{\partial \bq_i}\right)_i =
\end{equation*}

\footnotesize
$$  2\left( \begin {array}{c} {\it w_1}\,{\it v_1}\,{\it q_0}+{\it w_1}\,{\it v_2}\,{\it q_3}-{\it w_1}\,{\it v_3}\,{\it q_2}-{\it w_2}\,{\it v_1}\,{\it q_3}+{\it w_2}\,{\it v_2}\,{\it q_0}+{\it w_2}\,{\it v_3}\,{\it q_1}+{\it w_3}\,{\it v_1}\,{\it q_2}-{\it w_3}\,{\it v_2}\,{\it q_1}+{\it w_3}\,{\it v_3}\,{\it q_0}\\\noalign{\medskip}{\it w_1}\,{\it v_1}\,{\it q_1}+{\it w_1}\,{\it v_2}\,{\it q_2}+{\it w_1}\,{\it v_3}\,{\it q_3}+{\it w_2}\,{\it v_1}\,{\it q_2}-{\it w_2}\,{\it v_2}\,{\it q_1}+{\it w_2}\,{\it v_3}\,{\it q_0}+{\it w_3}\,{\it v_1}\,{\it q_3}-{\it w_3}\,{\it v_2}\,{\it q_0}-{\it w_3}\,{\it v_3}\,{\it q_1}\\\noalign{\medskip}-{\it w_1}\,{\it v_1}\,{\it q_2}+{\it w_1}\,{\it v_2}\,{\it q_1}-{\it w_1}\,{\it v_3}\,{\it q_0}+{\it w_2}\,{\it v_1}\,{\it q_1}+{\it w_2}\,{\it v_2}\,{\it q_2}+{\it w_2}\,{\it v_3}\,{\it q_3}+{\it w_3}\,{\it v_1}\,{\it q_0}+{\it w_3}\,{\it v_2}\,{\it q_3}-{\it w_3}\,{\it v_3}\,{\it q_2}\\\noalign{\medskip}-{\it w_1}\,{\it v_1}\,{\it q_3}+{\it w_1}\,{\it v_2}\,{\it q_0}+{\it w_1}\,{\it v_3}\,{\it q_1}-{\it w_2}\,{\it v_1}\,{\it q_0}-{\it w_2}\,{\it v_2}\,{\it q_3}+{\it w_2}\,{\it v_3}\,{\it q_2}+{\it w_3}\,{\it v_1}\,{\it q_1}+{\it w_3}\,{\it v_2}\,{\it q_2}+{\it w_3}\,{\it v_3}\,{\it q_3}\end {array} \right) 
 .$$ \normalsize \\

The vector obtained is quite ugly but one can see that it is linear in  $\bq$, it can thus be rewritten in a matrix-vector product:

\footnotesize
$$ 2  \underbrace{
\left( \begin {array}{cccc} {\it v_1}\,{\it w_1}+{\it v_2}\,{\it w_2}+{\it v_3}\,{\it w_3}&{\it v_3}\,{\it w_2}-{\it v_2}\,{\it w_3}&-{\it v_3}\,{\it w_1}+{\it v_1}\,{\it w_3}&{\it v_2}\,{\it w_1}-{\it v_1}\,{\it w_2}\\\noalign{\medskip}{\it v_3}\,{\it w_2}-{\it v_2}\,{\it w_3}&{\it v_1}\,{\it w_1}-{\it v_2}\,{\it w_2}-{\it v_3}\,{\it w_3}&{\it v_1}\,{\it w_2}+{\it v_2}\,{\it w_1}&{\it v_1}\,{\it w_3}+{\it v_3}\,{\it w_1}\\\noalign{\medskip}-{\it v_3}\,{\it w_1}+{\it v_1}\,{\it w_3}&{\it v_1}\,{\it w_2}+{\it v_2}\,{\it w_1}&{\it v_2}\,{\it w_2}-{\it v_1}\,{\it w_1}-{\it v_3}\,{\it w_3}&{\it v_2}\,{\it w_3}+{\it v_3}\,{\it w_2}\\\noalign{\medskip}{\it v_2}\,{\it w_1}-{\it v_1}\,{\it w_2}&{\it v_1}\,{\it w_3}+{\it v_3}\,{\it w_1}&{\it v_2}\,{\it w_3}+{\it v_3}\,{\it w_2}&{\it v_3}\,{\it w_3}-{\it v_1}\,{\it w_1}-{\it v_2}\,{\it w_2}\end {array} \right) 
}_{\mbox{\normalsize $\Delta[\v{v},\v{w}]$}}
 \left( \begin {array}{c} {\it q_0}\\\noalign{\medskip}{\it q_1}\\\noalign{\medskip}{\it q_2}\\\noalign{\medskip}{\it q_3}\end {array} \right) $$ \normalsize

By careful inspection of $\Delta[\v{v},\v{w}]$, we can identify a structure in the matrix that allows a compact notation

\begin{equation}
\Delta[\v{v},\v{w}] = 
\left( \begin {array}{cc}
\v{w} \cdot \v{v}    &      (\v{w} \c \v{v})^T   \\
\v{w} \c \v{v}         &       \v{w}\v{v}^T + \v{v}\v{w}^T  - \v{w} \cdot \v{v} \ I_3 
\end {array} \right).
\label{eqD_DELTA}
\end{equation}

That is

\begin{equation}
\frac{\partial (\v{v}^TR\v{w})}{\partial \bq} =
2 \Delta[\v{v},\v{w}] \bq
\label{eqD_DELTAq}
\end{equation}

And because $\v{v}^TR^T\v{w} = \v{w}^TR\v{v}$ we also have

\begin{equation}
\frac{\partial (\v{v}^TR^T\v{w})}{\partial \bq} =
2 \Delta[\v{w},\v{v}] \bq
\label{eqD_DELTAq2}
\end{equation}

\subsubsection{"Double $R$" Quadratic Form}

We are now interested in the derivative of a quadratic form involving $RJR^T$, that is, with the $\bq$ dependent matrix $R$ appearing twice. $J$ is an inertia matrix, therefore, $J=J^T$. This time, the vectors on the left ant on the right are the same, lets say $\v{u}$.

\begin{align*}
\half\frac{\partial }{\partial \bq} \left( \v{u}^TRJR^T\v{u}\right)  &=
\half \left(   \v{u}^T  \frac{\partial R}{\partial \bq_i}  JR^T   \v{u}   \right)_i  +
\half \left(   \v{u}^T RJ  \frac{\partial R^T}{\partial \bq_i}   \v{u}   \right)_i  \\
  &=   \left(   \v{u}^T  \frac{\partial R}{\partial \bq_i}  JR^T   \v{u}   \right)_i .
\end{align*}

Hence

\begin{equation}
\half\frac{\partial }{\partial \bq} \left( \v{u}^TRJR^T\v{u}\right)  = 2 \Delta[\v{u},JR^T\v{u}] \bq
\label{eqD_DELTAq3}
\end{equation}

\subsubsection{Properties}

By looking at \eqref{eqD_DELTA}, one may note the following relations

\begin{equation}
 \Delta[\v{v}_1+\v{v}_2,\v{w}] = \Delta[\v{v}_1,\v{w}] + \Delta[\v{v}_2,\v{w}] 
\label{eqD_DELTA_r1}
\end{equation}

\begin{equation}
 \Delta[\v{v},\v{w}_1+\v{w}_2] = \Delta[\v{v},\v{w}_1] + \Delta[\v{v},\v{w}_2] 
\label{eqD_DELTA_r2}
\end{equation}

\begin{equation}
 \Delta \left[\sum_{i=1}^{n}\v{v}_i,\sum_{j=1}^{m}\v{w}_j \right] = \sum_{i=1}^{n} \sum_{j=1}^{m} \Delta[\v{v}_i,\v{w}_j] 
\label{eqD_DELTA_r3}
\end{equation}

\begin{equation}
 \Delta[\alpha \v{v},\beta \v{w}] = \alpha \beta \Delta[\v{v},\v{w}]  
\label{eqD_DELTA_r4}
\end{equation}


\vspace{12pt}

\subsection{Time Derivative of $R$}

First note that by identification, one can verify that
\begin{equation}
G^T G = E^T E = I_4 - \b q \b q^T
\label{identity_in_R4}
\end{equation}
with $I_4$ the identity matrix in $\mathbb{R}^4$. Remember also 
\begin{equation*}
\Omega' = 2G\dot G^T = -2\dot G G^T \qquad \textrm{with} \qquad \Omega' \v v = \wp \times \v v
\end{equation*}
and
\begin{equation*}
\wp = 2G \dot{\b q} = -2\dot G \b q.
\end{equation*}

Now observe
\begin{align}
\Omega' R^T &= 2G\dot G^T G E^T \nonumber \\
                        &= -2\dot G G^T G E^T \nonumber \\
                        &=-2\dot G (I_4 -\b q \b q^T) E^T \nonumber \\
                        &=-2\dot G E^T -2\dot G \b q \underbrace{\b q^T E^T}_{(E\b q)^T=\v 0} \nonumber \\
                        &= -2\dot G E^T = -\dot R^T. \nonumber
\end{align}

We can finally write 

\begin{gather}
\dot R^T   = -\Omega' R^T      \label{dotRt_with_Omega} \\
\dot R       = -R \Omega'^T = R \Omega'.
\label{dotR_with_Omega}
\end{gather}



%
%
%


	

\pagebreak 

\graphicspath{{SubDocuments/SpeedComposition/}}

\section{Speed Composition}
\label{SpeedComposition}

Let be three referentials each designed by $0$, $1$ and $2$. Referential $0$ is inertial, referential $1$ is a rotating one and $2$ is a body fixed referential.\\
The same vector $\v{x}$ can be expressed in any of these referentials; when expressed in $0$, we will notate it as $\v{x}^0$, when expressed in $1$ it will be noted $\v{x}^1$ and $\v{x}^2$ in referential $2$. We will also write $x^i$ the quaternion $(0,\v{x}^i)$.\\
Moreover, three quaternions are defined: $q_{01}$ describes relative attitude of referential $1$ with respect to referential $0$,  $q_{12}$ describes relative attitude of referential $2$ with respect to referential $1$ and $q_{02}$  describes relative attitude of referential $2$ with respect to referential $0$.\\

\begin{center}
\psfrag{q01}[l][c][1][0]{$q_{01}$}
\psfrag{q12}[l][c][1][0]{$q_{12}$}
\psfrag{q02}[l][c][1][0]{$q_{02}$}
\psfrag{Ref0}[l][c][1][0]{$0$}
\psfrag{Ref1}[l][c][1][0]{$1$}
\psfrag{Ref2}[l][c][1][0]{$2$}
\epsfig{file=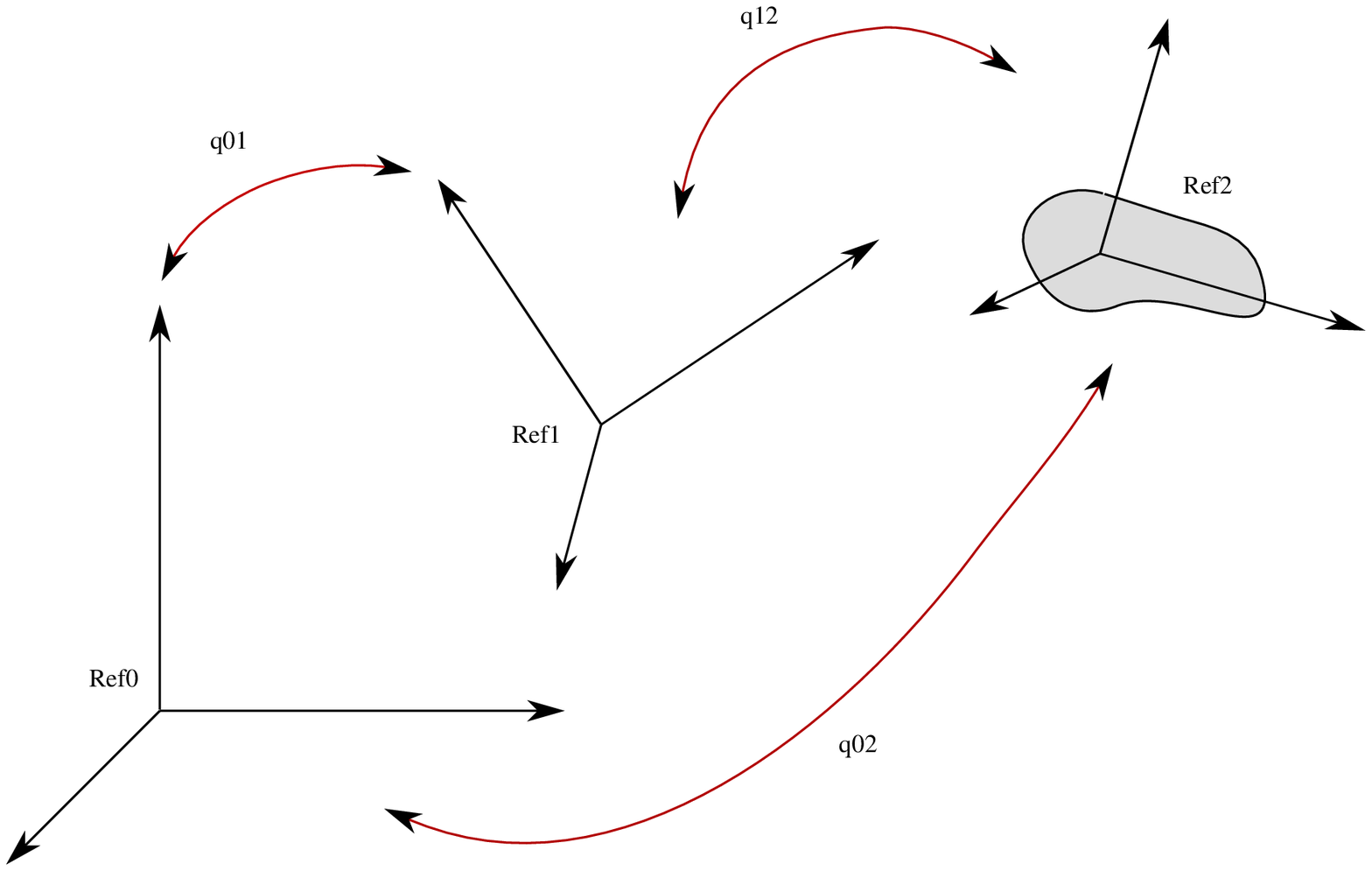, width=12cm}
\end{center}

So we may write

\begin{equation*}
x^0=q_{01} \o x^1 \o \bar{q}_{01}   \qquad   x^1=q_{12} \o x^2 \o \bar{q}_{12} \qquad    x^0=q_{02} \o x^2 \o \bar{q}_{02}
\end{equation*}

and by substitution

\begin{equation*}
x^0=q_{01} \o x^1 \o \bar{q}_{01}  =  q_{01} \o q_{12} \o x^2 \o \bar{q}_{12} \o \bar{q}_{01} 
  =  (q_{01} \o q_{12}) \o x^2 \o (\overline{q_{01} \o q_{12}})
\end{equation*}
  
we can identify $q_{02}$

\begin{equation}
q_{02} = q_{01} \o q_{12}.
\end{equation}

Noting $\omega_{ij}^j = (0,\v{\omega}_{ij}^j)$ the rotation velocity of the reference frame $j$ relative to frame $i$ expressed in the frame $j$ and remembering that $\omega_{ij}^j = 2\bar{q}_{ij} \o \dot{q}_{ij}$, we may write

\begin{align*}
\omega_{02}^2  &=   2\bar{q}_{02} \o \dot{q}_{02} \\
                              &=   2 (\bar{q}_{12} \o \bar{q}_{01}) \o (\dot{q}_{01} \o q_{12}    +   q_{01} \o \dot{q}_{12}) \\
                              &=   2 \bar{q}_{12} \o \bar{q}_{01} \o \dot{q}_{01} \o q_{12}    +
                                        2 \bar{q}_{12} \o \underbrace{\bar{q}_{01} \o  q_{01}}_{Id} \o \dot{q}_{12} \\
                              &=   \bar{q}_{12} \o \underbrace{(2 \bar{q}_{01} \o \dot{q}_{01})}_{\omega_{01}^1} \o q_{12}    +
                                         \underbrace{2 \bar{q}_{12} \o  \dot{q}_{12}}_{\omega_{12}^2}\\
                              &=  \bar{q}_{12} \o \omega_{01}^1 \o q_{12}    +      \omega_{12}^2 \\
                              &=   \omega_{01}^2   +      \omega_{12}^2.
\end{align*}

That is, we can add consecutive rotation speeds if they are expressed in the same referential.\\
In the case of the Cubsat, $\v{\omega}_{02}^2$ is the satellite's rotation velocity $\v{\omega}'$ expressed in body coordinates in the inertial referential model; we will note it $\v{\omega}'_{Inertial}$ here. On the other hand, $\v{\omega}_{12}^2$ is the satellite's rotation velocity $\v{\omega}'$ expressed in body coordinates in the non-inertial referential model (i.e. in orbital reference frame, ORF); we will note it $\v{\omega}'_{NonInertial}$. \\
$\v{\omega}_{01}^1$ is the ORF rotation velocity expressed in the ORF, that is $\v{\omega}_o$, while $\v{\omega}_{01}^2$ is the same vector, transformed in the body referential. This transformation is performed by $R^T$ from the non-inertial model ($\bar{q}_{12}$ in the above developement).\\
In other words, we can link the $\v{\omega}'$ vector from both inertial and non-inertial formulations (models) with

\begin{equation}
\v{\omega}'_{Inertial} = R^T_{NonInertial} \v{\omega}_o + \v{\omega}'_{NonInertial}.
\label{SC_comp}
\end{equation}

This is the speed to be used in computing the kinetic energy for the non-inertial model.


%
%
%


	

\pagebreak 

\graphicspath{{SubDocuments/Q2E_latex/}}



\section{Euler Angles to Quaternions}
\label{sec_Eul2Quat}

Three rotations by the Euler angles around each axis can be written as \\

$$ R_{\psi} =  \left[ \begin {array}{ccc} \cos \left( {\it \psi} \right) &-\sin \left( {\it \psi} \right) &0\\\noalign{\medskip}\sin \left( {\it \psi} \right) &\cos \left( {\it \psi} \right) &0\\\noalign{\medskip}0&0&1\end {array} \right] $$\\

$$ R_{\theta} =  \left[ \begin {array}{ccc} \cos \left( {\it \theta} \right) &0&\sin \left( {\it \theta} \right) \\\noalign{\medskip}0&1&0\\\noalign{\medskip}-\sin \left( {\it \theta} \right) &0&\cos \left( {\it \theta} \right) \end {array} \right] $$\\

$$ R_{\phi} =  \left[ \begin {array}{ccc} 1&0&0\\\noalign{\medskip}0&\cos \left( {\it \phi} \right) &-\sin \left( {\it \phi} \right) \\\noalign{\medskip}0&\sin \left( {\it \phi} \right) &\cos \left( {\it \phi} \right) \end {array} \right] $$\\

Combined together, they define the rotation matrix \\

$$ R = R_{\phi} R_{\theta} R_{\psi}. $$\\ 

Those three rotations can also be expressed as quaternion rotations \\

\begin{center}
$ \mathbf{q_{\phi}} =  \left[ \begin {array}{c} \cos \left( \frac{1}{2}\,{\it \phi} \right) \\\noalign{\medskip}\sin \left( \frac{1}{2}\,{\it \phi} \right) \\\noalign{\medskip}0\\\noalign{\medskip}0\end {array} \right] $ $\qquad$ $ \mathbf{q_{\theta}} =  \left[ \begin {array}{c} \cos \left( \frac{1}{2}\,{\it \theta} \right) \\\noalign{\medskip}0\\\noalign{\medskip}\sin \left( \frac{1}{2}\,{\it \theta} \right) \\\noalign{\medskip}0\end {array} \right] $ $\qquad$ $ \mathbf{q_{\psi}} =  \left[ \begin {array}{c} \cos \left( \frac{1}{2}\,{\it \psi} \right) \\\noalign{\medskip}0\\\noalign{\medskip}0\\\noalign{\medskip}\sin \left( \frac{1}{2}\,{\it \psi} \right) \end {array} \right] .$\\
\end{center}

\vspace{0.8cm}
The resulting quaternion can then be obtained by multiplying those three together \\

$$ \mathbf{q} = \mathbf{q_{\phi}} \circ \mathbf{q_{\theta}} \circ \mathbf{q_{\psi}} =  \left[ \begin {array}{c} \cos \left( \frac{1}{2}\,{\it \phi} \right) \cos \left( \frac{1}{2}\,{\it \theta} \right) \cos \left( \frac{1}{2}\,{\it \psi} \right) -\sin \left( \frac{1}{2}\,{\it \phi} \right) \sin \left( \frac{1}{2}\,{\it \theta} \right) \sin \left( \frac{1}{2}\,{\it \psi} \right) \\\noalign{\medskip}\cos \left( \frac{1}{2}\,{\it \psi} \right) \cos \left( \frac{1}{2}\,{\it \theta} \right) \sin \left( \frac{1}{2}\,{\it \phi} \right) +\cos \left( \frac{1}{2}\,{\it \phi} \right) \sin \left( \frac{1}{2}\,{\it \theta} \right) \sin \left( \frac{1}{2}\,{\it \psi} \right) \\\noalign{\medskip}\cos \left( \frac{1}{2}\,{\it \psi} \right) \cos \left( \frac{1}{2}\,{\it \phi} \right) \sin \left( \frac{1}{2}\,{\it \theta} \right) -\cos \left( \frac{1}{2}\,{\it \theta} \right) \sin \left( \frac{1}{2}\,{\it \phi} \right) \sin \left( \frac{1}{2}\,{\it \psi} \right) \\\noalign{\medskip}\cos \left( \frac{1}{2}\,{\it \phi} \right) \cos \left( \frac{1}{2}\,{\it \theta} \right) \sin \left( \frac{1}{2}\,{\it \psi} \right) +\cos \left( \frac{1}{2}\,{\it \psi} \right) \sin \left( \frac{1}{2}\,{\it \phi} \right) \sin \left( \frac{1}{2}\,{\it \theta} \right) \end {array} \right]. $$ \\

Note: that this result depends on the convention used in the order and choice of the Euler angles and rotation axes! \cite{mathworksQuat}


\pagebreak

\end{document}